\documentclass[12pt]{amsart}
\usepackage{booktabs} 
\usepackage[hidelinks]{hyperref}
\usepackage{enumerate}
\usepackage{enumitem}
\usepackage{xcolor}
\usepackage{float}
\allowdisplaybreaks
\usepackage{longtable}
\usepackage{amsfonts,amssymb,amsmath,textcomp}
\usepackage{float}
\usepackage{tikz-cd}
\usetikzlibrary{tikzmark}
\usepackage{tikz}
\usepackage{tabularx}
\usepackage{supertabular}
 2   

\usepackage{array,tabularx,xcolor} 

\newcolumntype{Y}{>{\raggedright\arraybackslash\color{blue}}X}

\usepackage[autostyle]{csquotes}
\theoremstyle{plain}
\newtheorem{thm}{Theorem}[section]

\newtheorem{lem}[thm]{Lemma}
\newtheorem{prop}[thm]{Proposition}
\newtheorem{cor}[thm]{Corollary}

\newtheorem{rem}[thm]{Remark}

\newcommand{\Gal}[1]{\text{Gal}(#1)}

\newcommand{\thmref}[1]{Theorem~\ref{#1}}

\theoremstyle{definition}
\newtheorem{defn}[thm]{Definition}

\newcommand{\Z}{\mathbb{Z}}
\newcommand{\Q}{\mathbb{Q}}
\newcommand{\R}{\mathbf{R}}
\newcommand{\C}{\mathbb{C}}

\begin{document}
\baselineskip 6mm

\title
{Absolutely Abelian Hilbert Class Fields and $\ell-$torsion conjecture}

\author{Mahesh Kumar Ram, Prem Prakash Pandey, Nimish Kumar Mahapatra}

\address[Mahesh Kumar Ram]{Indian Institute of Technology Madras, India.}
\email{maheshkumarram621@gmail.com}

\address[Prem Prakash Pandey]{Indian Institute of Science Education and Research Berhampur, India.}
\email{premp@iiserbpr.ac.in}

\address[Nimish Kumar Mahapatra]{Indian Institute of Technology Guwahati, India.}
\email{nimishravenshaw@gmail.com}


\date{\today}


\begin{abstract}
There are several recent works where authors have shown that number fields $K$ with `sufficiently many' units and cyclic class group contain a Euclidean ideal class provided the Hilbert class field $H(K)$ of $K$ is absolutely abelian. In this article, we explore the latter hypothesis: how often a number field $K$ has absolutely abelian Hilbert class field? For a number field $K$ to have absolutely abelian Hilbert class field, we obtain several criteria in terms of class number of $K$, P\'olya group of $K$, and genus number of $K$. We also show that for such number fields the $\ell-$torsion conjecture is true. Along with these, the article also reports some results on a theme to study class groups,  where primes of higher degree are used to study class groups.

\end{abstract}

\maketitle{}

\section{Introduction}

Class groups of number fields have been the subject of intensive research in the last two centuries. Class field theory, a gem of twentieth century mathematics, made it possible to view the class group of a number field $K$ as the Galois group of an extension $H(K)$, called the Hilbert class field of $K$, of the number field $K$. However, construction of Hilbert class fields of number fields is difficult. The mystery of class groups remains a subject of continuous investigation and several different approaches are developed to understand class groups.\\

In this article, we present some of our studies on class groups of number fields from several different approaches. Let $\mathcal{O}_K$ and $C\ell(K)$ denote the ring of integers and the class group of $K$, respectively. We use $D_K$ to denote the absolute value of the discriminant of $K$, and $n$ to denote the degree of $K$.  For any prime number $\ell$, we use $C\ell(K)[\ell]$ to denote the $\ell-$torsion subgroup of the class group $C\ell(K)$.  A very active and difficult theme \cite{BSTTTZ20, EV07, HV06, LB05, PCW20, WJ21, WJArxiv} is to obtain bounds of the form
\begin{equation}\label{ME}
|C\ell(K)[\ell]| \ll D_K^{\Delta+\epsilon},
\end{equation}
for arbitrary $\epsilon >0$ and some $0< \Delta <\frac{1}{2}$.  The constant in (\ref{ME}) depends only on $\ell, \epsilon$ and $n$ and not on $K.$ It is classical that one can take $\Delta=\frac{1}{2}$ (for example see \cite[Theorem 4.4]{WN80}). Motivated by the works of Brumer-Silverman \cite{BS96},  Duke \cite{WD98}, Zhang \cite{SZ05}, it is conjectured that one can take $\Delta=0$ \cite[Conjecture 7.1]{PCW20}. This conjecture is known as $\ell-$torsion conjecture and is connected to several counting problems (see \cite{BSTTTZ20, EV07, HV06, PCW20, PCW21}).  The survey article by Pierce \cite{LP22} is an excellent source to know more about the $\ell-$torsion conjecture.  From Gauss's genus theory, one can take $\Delta=0$ for the case $(n, \ell)=(2,2)$.  In a recent work, Kl\"{u}ners and Wang \cite{KW22} established the conjecture for $\ell-$torsion of class groups of $\ell-$extensions. For no other case the conjecture is known, in spite of tremendous attention given to the problem  \cite{BSTTTZ20, EPW17,  EV07, FW21, HV06,  KP22,  LTZ24, LB05, PCW20, AS17, WJ21}. From this illustrious list of work, we wish to highlight the important work of Venkatesh and Ellenberg \cite{EV07} which allows one to take $\Delta=\frac{1}{3}$ for the case $(n, \ell)=(2,3)$ and has been pivotal in many recent works.  We also would like to highlight the recent work of Chan-Koymans in this direction \cite{CK25}.\\

Before stating our result on $\ell-$torsion conjecture, we recall that the Hilbert class field $H(K)$ is called absolutely abelian if the extension $H(K)/\Q$ is abelian. Now we state one of the main results proved in this article.
\begin{thm}\label{Main}
Let $\ell$ be a prime number and $n$ be a positive integer. Every number field $K/\Q$ of degree $n$ whose Hilbert class field $H(K)$ is absolutely abelian satisfies 
\begin{equation}\label{ME1}
|C\ell(K)[\ell]| \ll_{\epsilon,\ell, n} D_K^{\epsilon}, 
\end{equation}
for every $\epsilon>0$. Here the constant depends only on $\ell, \epsilon$ and $n$,  and not on $K.$
\end{thm}
We prove Theorem \ref{Main} in Section 4.  In fact, we prove a stronger result, namely Theorem \ref{Main1}, from which Theorem \ref{Main} follows immediately.  It will be interesting to see how often number fields satisfy the hypothesis of Theorem \ref{Main}.  We are not aware of any literature explicitly discussing the following question:\\

\noindent {\it Question 1.} How often do number fields have absolutely abelian Hilbert class fields? \\

Question 1 is an important question by itself. But our interest in this question emanates from Theorem \ref{Main}, and an interesting series of recent works \cite{DGJ20, GM13, GS20,  KMJS22,SJ19} where a similar assumption is made. We briefly mention some of these below. If $\mathcal{O}_K$ is a Euclidean domain, then the class group $C\ell(K)$ is trivial. In 1972, under the generalized Riemann hypothesis, Weinberger \cite{Weinberger73} showed that the converse also holds for the number fields whose unit rank is at least one. Unconditionally, we know that for the number fields with unit rank at least four, the ring $\mathcal{O}_K$ is Euclidean if and only if the class group $C\ell(K)$ is trivial (see the work of Harper and Murty \cite{HM04}).  In 1979, Lenstra \cite{HWL79} introduced the notion of Euclidean ideal classes to study cyclicity of class groups of number fields. We urge the reader to read \cite{HWL79} for the definition of Euclidean ideal class.  Lenstra showed that existence of a Euclidean ideal class in $K$ ensures that $C\ell(K)$ is cyclic. Conversely, under generalized Riemann Hypothesis, Lenstra showed that number fields $K$ with unit rank at least one and with cyclic class group have a Euclidean ideal class. We state one unconditional result in this line of research due to Graves and Murty \cite{GM13}.

\begin{thm} \cite[Theorem 1]{GM13}\label{GM}
Let $K$ be a number field with unit rank at least $4$. Assume that the Hilbert class field $H(K)$ is abelian over $\Q$. If the class group $C\ell(K)$ is cyclic then it is generated by a Euclidean ideal class.
\end{thm}

There have been many recent works \cite{DGJ20, GS20,  KMJS22,SJ19} along the lines of the work of Graves and Murty. In all these works, the  authors assume that $K$ is a number field with absolutely abelian Hilbert class field $H(K)$. For such number fields, under some mild conditions (which vary from author to author) they show that if the class group $C\ell(K)$ is cyclic then $K$ has a  Euclidean ideal class. \\

The above discussion underlines the importance of Question 1. Our study on Question 1 is the subject matter of Section 2 and Section 3. In Section 2, we prove the following result.
\begin{thm}\label{OT}
Let $K$ be an abelian number field. Suppose $K$ contains a subfield $F$ with $[K : F] = n_1$. If the class group $C\ell(K)$ contains an element $\mathfrak{A}$ of order $m>1$ such that $m \nmid n_1 h_F$, then the extension $H(K)/\mathbb{Q}$ is not abelian.
\end{thm}

For quadratic number fields, we obtain a complete characterization of quadratic number fields with absolutely abelian Hilbert class fields. We prove the following result.
\begin{thm}\label{QN}
Let $K$ be a quadratic number field. Then $H(K)/\Q$ is abelian if and only if
$$C\ell(K) \cong (\Z/2\Z)^u,$$ for some $u \geq 0.$
\end{thm}

Another important result, relating to Question 1, we prove is Theorem \ref{T2}.  In Section 3, we study the roles of P\'olya groups  in relation to Question 1; theorems \ref{C1}, \ref{newrealquad}, \ref{realquad} are some results in this direction. \\

We remark that it is not known if there are infinitely many number fields $K$ with absolutely abelian Hilbert class field $H(K)$. On the other hand, let $H_{\ell}(K)$ denote the Hilbert $\ell-$class field of $K$, that is, $H_{\ell}(K)$ is a subfield of $H(K)$ such that the Galois group $\Gal{H_{\ell}(K)/K}$ is isomorphic to the $\ell-$part of $C\ell(K)$. Then it can be shown that $H_{\ell}(K)/\Q$ is abelian for infinitely many number fields $K$ (see Theorem \ref{infinite}). In Section 5 we present some results showing finiteness of some families of number fields with absolutely abelian Hilbert class fields. Based on Theorem \ref{QN}, we also list out real quadratic field with discriminant up to $10^5$ whose Hilbert class field is absolutely abelian.

%
%
%
Now we mention the next theme explored in this article. This is related to the following question:\\
{\it Question 2.} What are all pairs $(K,f)$ of number fields $K$ and positive integers $f$ such that the class group $C\ell(K)$ is generated by the ideal classes of primes of degree $f$?\\
For $f=1$, as a consequence of $\check{C}$ebotarev density theorem it follows that every number field $K$ can be taken. For $f>1$ the authors \cite{NPM24, PM22} have studied existence of such fields $K$, and illustrated that for such pairs $(K,f)$ one can obtain some further information on $C\ell(K)$ (see \cite{PM22, MR23, MR25}). Some of our recent findings on this theme are presented in Section 6.


{\bf Acknowledgement:} Initially we had obtained a weaker bound in (\ref{ME1}) as a consequence of Theorem \ref{T2}. We thank Jiuya Wang to point out that under our assumption we can prove the full $\ell-$torsion conjecture using Kronecker Weber Theorem.\\

\section*{Contents}
\phantomsection

\renewcommand{\tabularxcolumn}[1]{>{\color{blue}\raggedright\arraybackslash}m{#1}}

\renewcommand{\arraystretch}{1.2}

\begin{tabularx}{\textwidth}{X m{1.2cm}}
  2.~\hyperref[secA]{Absolutely abelian Hilbert class fields} & \pageref{secA} \\
  3.~\hyperref[secB]{P\'olya groups and genus fields} & \pageref{secB} \\
  4.~\hyperref[secC]{Proof of Theorem \ref{Main}} & \pageref{secC} \\
  
  5.~\hyperref[secE]{Some finiteness results and concluding remarks} & \pageref{secE} \\
  6.~\hyperref[secD]{Primes of higher degree} & \pageref{secD} \\
  7. ~ \hyperref[secF]{Appendix (Another proof of Theorem \ref{Main1})} & \pageref{secF}\\
  \hyperref[secG]{References} & \pageref{secG} \\
\end{tabularx}

\renewcommand{\arraystretch}{1}
\renewcommand{\tabularxcolumn}[1]{p{#1}}

\section{Absolutely abelian Hilbert class fields}\label{secA}
In this section we report our study relating to Question 1.  For any Galois extension $K/F$ of number fields, we use $G_{K/F}$ to denote the Galois group of the extension $K/F$. When $F=\Q$, we simply write $G_K$ instead of $G_{K/\Q}$. Then we have the following exact sequence
\begin{equation}\label{hpe1}
1 \longrightarrow G_{H(K)/K} \longrightarrow G_{H(K)/F} \longrightarrow G_{K/F} \longrightarrow 1.
\end{equation}
If $K/\Q$ is a cyclic extension of number fields then the exact sequence (\ref{hpe1}) always splits \cite{BW73}. Consequently, the Galois group $G_{H(K)}$ is a semi-direct product of $C\ell(K)$ and $G_K$. If this semi-direct product is a direct product then the Hilbert class field is absolutely abelian. This set up seems most obvious way to facilitate absolutely abelian Hilbert class fields. In this section we highlight families of number fields $K$ for which the Hilbert class field $H(K)$ can not be absolutely abelian or can be absolutely abelian. \\

We recall the following well known result. This can be proved in several ways (for example, see \cite{DP23, LW91}), for the sake of completeness we shall sketch a proof.
\begin{prop}\label{Con}
Let $K$ be a number field for which the Hilbert class field $H(K)$ is abelian. Then the conductors of $K$ and $H(K)$ are the same.
\end{prop}
We give a simple sketch based on a recent formula for the conductor of abelian number fields \cite{DP23}. Suppose degree of $K$ is 
$$n=2^{x_0}q_1^{x_1} \ldots q_s^{x_s},$$
where $x_0 \geq 0, x_1, \ldots, x_s>0$ and $q_i'$s are distinct primes. Let $p_1, \ldots, p_y$ be all the odd primes that ramify in $K$. If $2$ is unramified in $K$ then the conductor of $K$ is given by
\begin{equation}\label{e5}
p_1 \ldots p_y q_1^{\nu_{q_1}(e(q_1))} \ldots q_s^{\nu_{q_s}(e(q_s))}.
\end{equation}
Here $e(q_i)$ is ramification index of $q_i$ in $K$ and $\nu_q(z)$ for any prime $q$ and non-negative integer $z$ is the $q-$adic valuation of $z$. If $2$ ramifies in $K$ then the conductor of $K$ is 
\begin{equation}\label{e55}
2^{\nu_2(e(2))}p_1 \ldots p_y q_1^{\nu_{q_1}(e(q_1))} \ldots q_s^{\nu_{q_s}(e(q_s))} \mbox{ or } 2^{\nu_2(e(2))+1}p_1 \ldots p_y q_1^{\nu_{q_1}(e(q_1))} \ldots q_s^{\nu_{q_s}(e(q_s))}
\end{equation}
depending upon whether the ramification degree of $2$ in $K$ and $K(\sqrt{-1})$ are same or different.\\
As primes ramifying in $K$ and $H(K)$ are same and their ramification indices are also same, using (\ref{e5}) or (\ref{e55}), it is easy to see that both $K$ and $H(K)$ will have same conductor.\\

The following result can be deduced immediately from Proposition \ref{Con}.
\begin{cor}\label{cor1}
If $K$ is either a cyclotomic field or a real cyclotomic field or a subfield of $\Q(\zeta_{p^e})$. Then $H(K)/\Q$ is abelian if and only if $h_{K} = 1$.
\end{cor}
\begin{proof} Let $h_K = 1.$ Then, $K = H(K),$ and by assumption, the extension $K/\Q$ is abelian. Therefore, the extension $H(K)/\Q$ is also abelian. 

Now assume that $H(K)/\Q$ is abelian, and we will show that $h_K = 1.$ First, suppose that $K$ is a cyclotomic field, that is, $K = \Q(\zeta_m)$ for some positive integer $m \not \equiv 2 \pmod 4.$ Then, $m$ is the conductor of $K$. Proposition \ref{Con} implies that $m$ is also the conductor of $H(K).$ This forces $H(K) = K  = \Q(\zeta_m).$ Consequently,  $h_K = 1.$

Now assume that $K$ is a real cyclotomic field with conductor $m$, that is, $K = \Q(\zeta_m + \zeta_m^{-1})$. By Proposition \ref{Con}, the conductor of $H(K)$ is also $m$. Therefore, $H(K) \subset \Q(\zeta_m)$. Using this, together with the facts that $K \subset H(K)$ and $[\Q(\zeta_m) : K] = 2$, it follows that either $H(K) = K$ or $H(K) = \Q(\zeta_m)$.  Since the extension $H(K)/K$ is unramified at all infinite primes and $K \subset \mathbb{R}$, we must have $H(K) = K$. Consequently, $h_K = 1$.

Now, assume that $K$ is a subfield of $\Q(\zeta_{p^e})$. Then the conductor of $K$ is $p^u$, where $0 \leq u \leq e$. Proposition \ref{Con} ensures that the conductor of $H(K)$ is also $p^u$. Therefore, $$K \subset H(K) \subset \Q(\zeta_{p^u}).$$ Since $p$ is unramified in the extension $H(K)/K$ but $p$ is totally ramified in the extension $\Q(\zeta_{p^u})/\Q$, it follows that $K = H(K)$. Hence $h_K = 1$, which completes the proof.
\end{proof}
The following result is one of the main result we prove in this section.
\begin{thm}\label{T2}
Let $S$ be any finite set of primes. Let $K$ be an abelian number field whose conductor has prime factors only from the set $S$.  Assume that $G_{H(K)}$ is abelian, then there exists a number $t$ depending only on $S$ such that the following holds\\
(a) all the prime factors of $h_K$ are factors of $t$,\\
(b) the class number of $K$ satisfies the bound $h_K \leq t$, and  moreover $h_K \mid t$.

\end{thm}
\begin{rem}
Having infinitely many number fields whose conductors are divisible only by primes in a finite set $S$ of primes and whose Hilbert class fields are absolutely abelian facilitates the existence of infinitely many number fields with a fixed class number.
\end{rem}

Before proving Theorem \ref{T2}, we introduce some notations that will be used in the sequel. For an extension $K/F$ of number fields, if $\mathfrak{P}$ is a prime ideal of $K$ lying above a prime ideal $\mathfrak{p}$ of $F$, then we use $e(\mathfrak{P}|\mathfrak{p})$ to denote the ramification index of $\mathfrak{P}$ for the extension $K/F$. When $K/F$ is Galois, then $e(\mathfrak{P}|\mathfrak{p})$ does not depend on the prime $\mathfrak{P}$. Thus, when we do not need to keep track of the prime $\mathfrak{P}$, we should simply use $e(\mathfrak{p})$ instead of $e(\mathfrak{P}|\mathfrak{p})$. Now we are in a position to prove Theorem \ref{T2}
\begin{proof}[Proof of Theorem $\ref{T2}$]
Suppose $S=\{p_1, \ldots, p_r\}$ and let $S_1=\{q_1, \ldots, q_s\}$ denote the exact set of primes which divide at least two elements from the set $\{p_1-1, \ldots, p_r-1\}$. Let $u_i$ be the largest integer such that $q_i^{u_i}|\prod_j (p_j-1)$. Now suppose $m=p_1^{v_1} \ldots p_r^{v_r}$ is the conductor of $K$ for non-negative integers $v_1, \ldots, v_r$. We write $x$ for the product of $p_i$ for which $v_i \geq 2$ and $p_i$ divides exactly one of $p_1-1, \ldots, p_r-1$.  Next, let $w_i$ be the largest integer such that $p_i^{w_i}|\prod_j (p_j-1)$. We define 
$$t=\prod_i q_i^{u_i} \prod_{p_i|x}p_i^{w_i}.$$
We remark that $t=1$ when $S_1=\emptyset$ and $x=1$. Note that $t|\prod_j(p_j-1)$ and there may exist an $i \in\{1, \ldots, r\}$ such that $p_i \nmid x$ but $p_i|t$.
As $H(K)/\Q$ is abelian, from Proposition \ref{Con}, we see that
\begin{equation}\label{e6}
\Q \subset K \subset H(K) \subset \Q(\zeta_m).
\end{equation}
First, we prove (a). Let $\ell$ be a prime divisor of $h_K$ that does not divide $t$. As $h_K=[H(K):K],$ from (\ref{e6}) it follows that $\ell | [\Q(\zeta_m):\Q]$. That is,
\begin{equation}\label{e7}
\ell |\phi(m)=\prod_i p_i^{(v_i-1)}(p_i-1).
\end{equation}
Now we consider three cases.\\
Case- (i): $\ell =p_i$ for some $i$ with $v_i \geq 2$.\\
Since $\ell$ does not divide $t$, we conclude that $\ell$ does not divide $p_j-1$ for any $j$. Thus we conclude that the largest power of $p_i$ dividing $[\Q(\zeta_m):\Q]$ is $p_i^{(v_i-1)}$. Further, $\ell=p_i$ divides $[H(K):K]$, hence the largest power of $p_i$ dividing $[K:\Q][\Q(\zeta_m):H(K)]$ is strictly smaller than $p_i^{(v_i-1)}$.\\
Suppose $\wp_i$ is a prime ideal of $\Q(\zeta_m)$ lying above $p_i$, and let $\mathfrak{P}_i, \mathfrak{p}_i$ denote the primes of $H(K)$ and $K$ respectively below $\wp_i$. It is well known that
\begin{equation}\label{e8}
e(\wp_i|p_i)=p_i^{(v_i-1)}(p_i-1).
\end{equation}
Also, we have 
\begin{equation}\label{e9}
e(\wp_i|p_i)=e(\wp_i|\mathfrak{P}_i)e(\mathfrak{P}_i|\mathfrak{p}_i)e(\mathfrak{p}_i|p_i).
\end{equation}
Note that $e(\mathfrak{P}_i|\mathfrak{p}_i)=1$,  $e(\wp_i|\mathfrak{P}_i)| [\Q(\zeta_m):H(K)]$ and $e(\mathfrak{p}_i|p_i)|[K:\Q]$. Hence we conclude that the power of $p_i$ in $e(\wp_i|p_i)$ in (\ref{e8}) and (\ref{e9}) do not match. This contradiction establishes that there is no prime $\ell$ dividing $h_K$ and not dividing $t$.\\
Case- (ii): $\ell = p_i$ for some $i$ with $v_i=1$.\\ 
From (\ref{e7}), we must have $\ell |(p_j-1)$ for some $j$. Since $\ell$ does not divide $t$, there exist a unique $j$ such that $\ell|(p_j-1)$. Suppose that the largest power of $\ell$ dividing $(p_j-1)$ is $\ell^c$ for some positive integer $c$. Then the largest power of $\ell$ dividing $[\Q(\zeta_m):\Q]$ is also $\ell^c$. Further, $\ell$ divides $[H(K):K]$, hence the largest power of $\ell$ dividing $[K:\Q][\Q(\zeta_m):H(K)]$ is strictly smaller than $\ell^c$.\\
As earlier, let $\wp_j$ be a prime ideal of $\Q(\zeta_m)$ lying above $p_j$, and let $\mathfrak{P}_j, \mathfrak{p}_j$ denote the primes of $H(K)$ and $K$ respectively below $\wp_j$. Then we have
\begin{equation}\label{e10}
e(\wp_j|p_j)=p_j^{(v_j-1)}(p_j-1).
\end{equation}
Also, we have 
\begin{equation}\label{e11}
e(\wp_j|p_j)=e(\wp_j|\mathfrak{P}_j)e(\mathfrak{P}_j|\mathfrak{p}_j)e(\mathfrak{p}_j|p_j).
\end{equation}
Note that $e(\mathfrak{P}_j|\mathfrak{p}_j)=1$,  $e(\wp_j|\mathfrak{P}_j)| [\Q(\zeta_m):H(K)]$ and $e(\mathfrak{p}_j|p_j)|[K:\Q]$. Hence we conclude that the power of $\ell$ dividing $e(\wp_j|p_j)$ in (\ref{e10}) is different from the power of $\ell$ dividing $e(\wp_j|p_j)$ in (\ref{e11}). This contradiction establishes that there is no prime $\ell$ dividing $h_K$ and not dividing $t$.\\
Case- (iii): $\ell \neq p_i$ for any $i$.\\
Once again we see that $\ell$ divides $p_j-1$ for precisely one $j$. Now proof goes exactly similar as in case- (ii).\\
Now we prove (b).\\
From the proof of (a) we know that all the prime factors of $h_K$ are factors of $t$. Let $\ell$ be any prime dividing $h_K$. Then $\ell | \prod_i p_i^{(v_i-1)} (p_i-1)$. Once again we consider different cases.\\
Case- (i): $\ell =p_i$ for some $i$ with $v_i \geq 2$.\\
Suppose $\ell^a$ is the largest power of $\ell$ dividing $h_K$ and $\ell^b$ is the largest power of $\ell$ dividing $t$. As $t|\phi(m)$, we see that the largest power of $p_i$ dividing $[\Q(\zeta_m):\Q]$ is $p_i^{(b+v_i-1)}$. Further, the largest power of $\ell=p_i$ dividing $[H(K):K]$ is $\ell^a$, hence the largest power of $p_i$ dividing $[K:\Q][\Q(\zeta_m):H(K)]$ is $p_i^{(b+v_i-1-a)}$.\\
Suppose $\wp_i$ is a prime ideal of $\Q(\zeta_m)$ lying above $p_i$, and let $\mathfrak{P}_i, \mathfrak{p}_i$ denote the primes of $H(K)$ and $K$ respectively below $\wp_i$. It is well known that the ramification index
\begin{equation}\label{e12}
e(\wp_i|p_i)=p_i^{(v_i-1)}(p_i-1).
\end{equation}
Also, we have 
\begin{equation}\label{e13}
e(\wp_i|p_i)=e(\wp_i|\mathfrak{P}_i)e(\mathfrak{P}_i|\mathfrak{p}_i)e(\mathfrak{p}_i|p_i).
\end{equation}
Note that $e(\mathfrak{P}_i|\mathfrak{p}_i)=1$,  $e(\wp_i|\mathfrak{P}_i)| [\Q(\zeta_m):H(K)]$ and $e(\mathfrak{p}_i|p_i)|[K:\Q]$. Hence the largest power of $p_i$ in $e(\wp_i|p_i)$ is at most $p_i^{(b+v_i-1-a)}$. Comparing this with (\ref{e12}) we conclude that $b \geq a$. \\
Case- (ii): $\ell = p_i$ for some $i$ with $v_i=1$ or $\ell \neq p_i$ for any $i$.\\
Since $\ell |t$, we must have $\ell$ divides $p_j-1$ for more than one $j$.  We see that the largest power of $\ell$ dividing $t$ and the largest power of $\ell$ dividing $\phi(m)=[\Q(\zeta_m):\Q]$ are same, that is, $\ell^b$. But $h_K| [\Q(\zeta_m):\Q]$, and hence the highest power of $\ell$ dividing $[\Q(\zeta_m):\Q]$ is greater than ore equal to $\ell^a$. This establishes that $h_K \leq t$.
\end{proof}

We can deduce several corollaries, analogous to Corollary \ref{cor1}, from Theorem \ref{T2}. 
\begin{cor}\label{cor4}
Let $m=\prod_{i=1}^r p_i^{\nu_i}$, where $p_i'$s are distinct primes such that no $p_i$ divides $p_j-1$ and the only common factor between $p_i-1$ and $p_j-1$ is a power of $2$.  Suppose $K$ is a subfield of $\Q(\zeta_m)$ of odd degree. Then, the extension $H(K)/\Q$ is abelian if and only if $h_K = 1$.
\end{cor}

\begin{proof} Assume that $H(K)/\Q$ is abelian.  Let $S=\{p_1, \cdots, p_r\}$, and define the set $S_1$ and $t$ as in the proof of Theorem \ref{T2}. It is readily seen that $t=2^s$ for some non-negative integer $s$. From Theorem \ref{T2}, we obtain
$$h_K=2^u \mbox{ for some }u \leq s.$$ As the degree $[K:\Q]$ is odd, from Corollary \ref{AC}, we see that $2 \nmid h_K$. Thus $h_K=1.$

Conversely, assume that $h_K = 1$. Then $K = H(K)$. Since $K$ is assumed to be a subfield of $\mathbb{Q}(\zeta_m)$, it follows that $H(K)/\mathbb{Q}$ is an abelian extension. This completes the proof.
\end{proof}
In a similar manner, we obtain the following.
\begin{cor}\label{cor5}
Let $K$ be an abelian number field of conductor $m$. Suppose $\phi(m)=[K:\Q].d,$ where $(d, [K:\Q])=1.$ Then $H(K)/\Q$ is abelian if and only if $h_K=1.$
\end{cor}
\begin{proof}
Assume that $H(K)/\Q$ is abelian. By Proposition \ref{Con}, we have $H(K) \subset \Q(\zeta_m)$. Consequently,
$$[\Q(\zeta_m):H(K)][H(K):K][K:\Q]=\phi(m)=[K:\Q].d.$$
This shows that $[H(K):K] |d.$ An immediate application of Corollary \ref{AC} implies that $h_K=1.$
The converse is straightforward.
\end{proof}

Before proceeding with the proof of Theorem \ref{OT}, we introduce some additional notations and recall some results from class field theory.  For any unramified prime ideal $\mathfrak{P}$ of $K$, the unique element $\sigma$ of the decomposition group $D_{\mathfrak{P}} \subset G_{K/F}$ satisfying 
$$\sigma (x) \equiv x^{N(\mathfrak{p})} \pmod {\mathfrak{P}}\; \; \forall x \in \mathcal{O}_K$$
is called the Frobenius automorphism of $\mathfrak{P}$ for the extension $K/F$ and is denoted by $\left(\frac{\mathfrak{P}}{K/F}\right)$, and by $\left(\frac{\mathfrak{p}}{K/F}\right)$ in case $G_{K/F}$ is abelian.  Here $\mathfrak{p}=\mathfrak{P} \cap \mathbf{O}_F$ and $N(\mathfrak{p})$ is the norm for the extension $F/\Q$. The residue degree of the prime $\mathfrak{P}$ for the extension $K/F$ is denoted by $f(\mathfrak{P}|\mathfrak{p})$. We recall the following result from \cite{GJ96}.
\begin{lem}\cite[p.126]{GJ96}\label{FR}
Let $K/F$ be a Galois extension of number fields and $E$ be an intermediate field. For any unramified prime ideal $\mathfrak{P}$ of $K$, we have
$$\left( \frac{\mathfrak{P}}{K/F} \right)^{f(\mathbf{p}/\mathfrak{p})}=\left( \frac{\mathfrak{P}}{K/E}\right).$$
Here $\mathbf{p}$ and $\mathfrak{p}$ are the prime ideals below $\mathfrak{P}$ in $E$ and $F$ respectively.
\end{lem}
A very important result of the class field theory is the following \cite{GJ96}.
\begin{thm}\cite[Theorem 12.1]{GJ96}\label{HCF}
The Galois group $G_{H(K)/K}$ is isomorphic to the class group $C\ell(K)$. The isomorphism is induced by the Artin map
$$[\mathfrak{P}] \longmapsto \left( \frac{\mathfrak{P}}{H(K)/K} \right).$$
\end{thm}

Next, we recall the $\check{C}$ebotarev density Theorem. For any $\sigma \in G_{K/F}$ we consider the set $P_{K/F}(\sigma)$ of prime ideals $\mathfrak{p}$ of $F$ such that there is a prime $\mathfrak{P}$ of $K$ above $\mathfrak{p}$ satisfying
$$\sigma =\left( \frac{\mathfrak{P}}{K/F} \right),$$
also we use $C_{\sigma}$ to denote the conjugacy class of $\sigma$ in $G_{K/F}$. Now we state the $\check{C}$ebotarev density theorem.
\begin{thm}\cite[Theorem 5.2]{GJ96}\label{CDT}
Let $K/F$ be a Galois extension with Galois group $G_{K/F}$. For every $\sigma \in G_{K/F}$, the density of the set $P_{K/F}(\sigma)$ is positive and equals to $\frac{|C_{\sigma}|}{[K:F]}$. 
\end{thm}
The following result can be derived from Theorem \ref{HCF} and Theorem \ref{CDT}. A completely algebraic proof of slightly weaker statement can be found in \cite{LS91}.
\begin{thm}\label{CCT} Each ideal class of a number field $K$ contains infinitely many prime ideals of $K$ of residue degree $1$.
\end{thm}
\begin{proof}[Proof of Theorem $\ref{OT}$]
  It suffices to show that $H(K)/F$ is non-abelian. Suppose, for the sake of contradiction,  that $H(K)/F$ is abelian.  Theorem \ref{CCT} implies that $\mathfrak{A}$ contains a prime ideal $\mathfrak{P}$ of $K$ of residue degree one, such that $\mathfrak{P}$ is unramified in $K/\Q$ . Put $\mathfrak{p} = \mathfrak{P}\cap \mathcal{O}_F$. Clearly, $\mathfrak{p}$ is not inert in $K$. By Dedekind's theorem on the prime factorization of an ideal, we have
\begin{equation}\label{Me}
    \mathfrak{p} \mathcal{O}_K = \mathfrak{P}_1\mathfrak{P}_2\cdots \mathfrak{P}_{n_1}
\end{equation}
 where $\mathfrak{P} =\mathfrak{P}_1$ and the primes $\mathfrak{P}_i$ are distinct. Let 

\begin{equation}\label{Me1}
\mathfrak{P}_i \mathcal{O}_{H(K)} = \prod_{j=1}^{r_i} \mathfrak{P}_{ij} \; , 
\end{equation}
where $\mathfrak{P}_{ij}$ are prime ideals of $H(K)$ lying above $\mathfrak{P}_i$. Consequently,

$$\mathfrak{p} \mathcal{O}_{H(K)} = \prod_{i=1}^{n_1}\prod_{j=1}^{r_i} \mathfrak{P}_{ij}.$$ This shows that the prime ideals $\mathfrak{P}_{ij}$ are all the primes of $H(K)$ lying above $\mathfrak{p}$. Since $H(K)/F$ is abelian, it follows that
\begin{equation}\label{Me2}
\left(\frac{\mathfrak{P}_{ij}}{H(K)/F} \right) 
\end{equation} 
are all same for all $i$ and $j$, and can be denote by $\left(\frac{\mathfrak{p}}{H(K)/F} \right)$. Using Lemma \ref{FR}, for each fixed $i$, we have 

\begin{equation}\label{Me3}
	\left(\frac{\mathfrak{P}_{ij}}{H(K)/F} \right)^{f(\mathfrak{P}_i/{\mathfrak{p})}} = \left(\frac{\mathfrak{P}_{ij}}{H(K)/K} \right)
	\end{equation} 
Since $H(K)/K$ is abelian, for each fixed $i$ we have

\begin{equation}\label{Me4}
	\left(\frac{\mathfrak{P}_{ij}}{H(K)/K} \right) = \left(\frac{\mathfrak{P}_{i}}{H(K)/K} \right)
\end{equation} 
for all $j$.
Using identities \eqref{Me2}- \eqref{Me4}, and the fact
$$f(\mathfrak{P}_1/{\mathfrak{p}})= f(\mathfrak{P}_2/{\mathfrak{p}})= \cdots = f(\mathfrak{P}_{n_1}|{\mathfrak{p}}) = 1,$$  we have
$$\left(\frac{\mathfrak{P}_{1}}{H(K)/K} \right)= \left(\frac{\mathfrak{P}_{2}}{H(K)/K} \right)= \cdots = \left(\frac{\mathfrak{P}_{n_1}}{H(K)/K} \right) = \tau, (\text{say}).$$ \\
By Theorem \ref{HCF}, the order of $\tau$ is equal to the order of $\mathfrak{A}$. Therefore, the order of $\tau$ is $m$. Also,  $\mathfrak{P_1}^m$ is a principal ideal. Consequently,
\begin{align*}
N_{K/F}(\mathfrak{P_1}^m) &= (N_{K/F}(\mathfrak{P_1}))^m \\ &= \left(\left(\prod_{\sigma \in G_{K/F}} \sigma(\mathfrak{P}_1)\right) \cap \mathcal{O}_F\right)^m \\  & =\mathfrak{p}^m
\end{align*}
is a principal ideal in $\mathcal{O}_F$. If $u = \text{gcd}(m, h_F)$, then $\mathfrak{p}^u$ is a principal ideal. From \eqref{Me}, we conclude  that $\left(\prod_{i =1}^{n_1} \mathfrak{P}_i\right)^u$ is a principal ideal in $\mathcal{O}_K$. Applying Theorem \ref{HCF}, we obtain 
$$\left(\prod_{i=1}^{n_1} \left(\frac{\mathfrak{P}_{i}}{H(K)/K} \right) \right)^u = \tau^{n_1u} = Id.$$ This shows that the order of $\tau$ divides $n_1u$, which contradicts the assumption that $m $ does not divide $n_1h_F$. This shows that the extension $H(K)/\Q$ is not abelian. \\

\end{proof}

As an immediate corollary to Theorem \ref{OT}, we have the following.
\begin{cor}\label{AC}
Let $K$  be an abelian number field of degree $n$. If there exists a prime $\ell$ such that $\ell$ divides $h_K$ but does not divide $n$, then the extension $H(K)/\Q$ is not abelian.
\end{cor}


The following result is another consequence of Theorem \ref{OT}.
\begin{cor} \label{Th23} Let $n>1$  and $m>1$ be integers such that there is a prime factor $p$ of $m$ which does not divide $n$. Suppose $K/\Q$ is a cyclic extension of number fields of degree $n$, and $\text{Aut}(G)$ denotes the automorphism group of a group $G$. If $n$ is relatively prime to $|\text{Aut}(G)|$ for each abelian group $G$ of order $m$, then $h_K \neq m$.
\end{cor} 
\begin{proof} Assume, by contradiction, that $h_K = m.$ Since $K/\Q$ is a cyclic extension, as mentioned at the beginning of this section, the sequence
\begin{equation*}
1 \longrightarrow C\ell(K) \longrightarrow G_{H(K)} \longrightarrow G_{K} \longrightarrow 1
\end{equation*} splits. Thus, the group $G_{H(K)}$ must be a semidirect product $C\ell(K) \rtimes_{\phi} G_{K}$ for some homomorphism $\phi: G_K \longrightarrow Aut(C\ell(K))$. By the assumption that $n$ is relatively prime to $|\text{Aut}(G)|$ for each abelian group $G$ of order $m$, it follows that $\phi$ is a trivial group homomorphism. Thus, 
\begin{equation}\label{E3}
    G_{H(K)} \simeq C\ell(K) \oplus G_K.
\end{equation}
Since the groups $ C\ell(K)$ and  $G_K$ are abelian, the identity \eqref{E3} ensures that $G_{H(K)}$ is also abelian, that is, the extension $H(K)/\Q$ is abelian. This contradicts Corollary \ref{AC}, because we assumed that $m$ has a prime factor $p$ that does not divide $[K:\Q]$. Therefore, the class number of $K$ cannot be $m$. This completes the proof.
\end{proof}

We mention one more consequence of Theorem \ref{OT}.
\begin{cor}\label{cor24}
Let $p$  and $q$ be two distinct primes such that $p$ and $q-1$ are relatively prime. Suppose $e$ is a positive integer. If $K/\Q$ is a cyclic extension of degree $p^e$, then $h_K \neq q$.
\end{cor}
From Corollary \ref{cor24}, it follows readily that if $K/\Q$ is a cyclic extension of prime-power degree, then the class number of $K$ cannot be $2$. \\

The $2$-part of class groups of quadratic fields is well investigated \cite{YTC20}, and are related to the solvability of negative Pell equations.  From genus theory, there is a good understanding of the $2$-rank of class groups. As a consequence of Corollary \ref{AC}, we deduce the following result.
\begin{cor}\label{C2}
Let $K$ be a quadratic number field and $m$ be its conductor. Suppose $2^d$ is the highest power of $2$ dividing $\phi(m)$. If $H(K)$ is absolutely abelian, then the class number $h_K=2^s$ for some integer $0 \leq s \leq d-1$.
\end{cor}

\begin{proof}
From Corollary \ref{AC} it follows that 
\begin{equation*}
|C\ell(K)|=2^s \mbox{ for some non-negative integer }s.
\end{equation*}
Further we see that $\Q \subset K \subset H(K) \subset \Q(\zeta_m)$. Thus, the highest power of $2$ dividing $[H(K):K]$ is at most $2^{d-1}$. Now the result follows.
\end{proof}

Now we prove Theorem \ref{QN}.
\begin{proof}[Proof of Theorem \ref{QN}]
Assume that $H(K)/\Q$ is an abelian extension. By Theorem \ref{OT}, the order of every element of $C\ell(K)$ divides $2$. Since $C\ell(K)$ is abelian, it follows that
$
C\ell(K) \simeq (\Z/2\Z)^u
$
for some integer $u \geq 0$.\\

Conversely, assume that $C\ell(K) \simeq (\Z/2\Z)^u$
for some  $u \geq 0$. Since $K$ is a quadratic field,
\[
G_K=\langle \tau\rangle \simeq \Z/2\Z,
\]
where $\tau$ denotes the non-trivial automorphism of $K$ over $\Q$. Let $\tilde{\tau}\in G_{H(K)}$ be a lift of $\tau$. Let $\mathfrak{p}$ be an unramified prime ideal of $\mathcal{O}_K$, and let
\[
\sigma_{\mathfrak{p}}
=
\left(\frac{\mathfrak{p}}{H(K)/K}\right)\in G_{H(K)/K}
\]
denote its Frobenius automorphism. By the properties of the Artin symbol, we have
\begin{equation}\label{Nie20}
    \tilde{\tau}\sigma_{\mathfrak{p}}\tilde{\tau}^{-1}
=
\left(\frac{\tau(\mathfrak{p})} {H(K)/K}\right).
\end{equation}
Since $\mathfrak{p}$ and $\tau(\mathfrak{p})$ are conjugate over $\Q$, the ideal $\mathfrak{p}\tau(\mathfrak{p})$ is a principal ideal of $\mathcal{O}_K$. Since the Artin map is trivial on principal ideals, we obtain
\[
\left(\frac{\mathfrak{p}\tau(\mathfrak{p})}{H(K)/K}\right)=Id.
\]
Therefore,
\begin{equation}\label{Nie21}
    \left(\frac{\tau(\mathfrak{p})}{H(K)/K}\right)
=
\left(\frac{\mathfrak{p}}{H(K)/K}\right)^{-1}
=
\sigma_{\mathfrak{p}}^{-1}.
\end{equation}
From \eqref{Nie20} and \eqref{Nie21}, we have
\[
\tilde{\tau}\sigma_{\mathfrak{p}}\tilde{\tau}^{-1}
=
\sigma_{\mathfrak{p}}^{-1}.
\]
Since the Frobenius elements $\sigma_{\mathfrak{p}}$ generate $G_{H(K)/K}$, the above inversion relation holds for every element $\sigma\in G_{H(K)/K}$. That is, 
\begin{equation}\label{Nie22}
    \tilde{\tau}\sigma\tilde{\tau}^{-1}
=
\sigma^{-1}
\qquad
\text{for all } \sigma\in G_{H(K)/K}.
\end{equation}
By hypothesis,
\[
G_{H(K)/K}\simeq C\ell(K)\simeq (\Z/2\Z)^u.
\]
Thus, every element $\sigma\in G_{H(K)/K}$ satisfies $\sigma^2=1,$
which is equivalent to $\sigma=\sigma^{-1}$. Consequently, \eqref{Nie22} becomes
\[
\tilde{\tau}\sigma\tilde{\tau}^{-1}
=
\sigma
\]
for all $\sigma\in G_{H(K)/K}$. Hence, $\tilde{\tau}$ commutes with every element of $G_{H(K)/K}$. Thus, using the fact that  $[G_{H(K)}:G_{H(K)/K}]=2$ and $G_{H(K)/K}$ is abelian, we conclude that  $G_{H(K)}$ is abelian. Consequently, $H(K)/\Q$ is an abelian extension. This completes the proof.
\end{proof}

\section{P\'olya groups and genus fields}\label{secB}
Let $Po(K)$ denote the P\'olya group of $K$.  P\'olya groups are a subgroup of class group and are easier to describe, for more details on P\'olya groups we refer to \cite{BR63, PJC, JLC, JLC1, SRM22, ZAN}.  It looks like that, for number fields $K$ with absolutely abelian Hilbert class field $H(K)$, the P\'olya group $Po(K)$ of $K$ covers the class group $C\ell(K)$ pretty well. In this direction, we prove the following result.
\begin{thm}\label{C1}
Let $K$ be a cyclic number field of degree $n$. Assume that either $n$ is odd or $K$ is not contained in $\mathbb{R}$. Then, the extension $H(K)/\Q$ is abelian if and only if $Po(K)=C\ell(K)$.
\end{thm}
We prove Theorem \ref{C1} using a result of Chabert \cite{JLC} and some results on genus fields. There are several studies on genus fields of number fields \cite{MI74, MI75, ZX85, FY67}. For Ishida \cite{MI74} and Zhang \cite{ZX85}, the genus field $K^*$ of $K$ is the maximal abelian extension of $K$ which is a composite of $K$ with an abelian extension $K_0^{*}$ of $\Q$ and is unramified at all the finite primes. Whereas, for Furuta \cite{FY67}, the genus field $\widetilde{K}$ of $K$ is the maximal abelian extension of $K$ which is composite of $K$ with an abelian extension $\tilde{K_0}$ of $\Q$ and is unramified at all primes. In this article the former genus field will be referred as the narrow genus field and the latter will be referred as the absolute genus field. It is known that the degrees $[K^*:K]$ and $[\tilde{K}:K]$ differ by a power of $2$ (for example see page 2216 in \cite{FLN23}). However, we are not aware of any literature classifying all number fields for which the narrow genus field and absolute genus field are same. We show that for a cyclic number field $K$ of odd degree,  or an imaginary cyclic number field $K$, the narrow genus field and the absolute genus field coincide. (see the proof of Theorem \ref{C1}). \\

Before proving Theorem \ref{C1}, we recall a result of Chabert \cite{JLC}.
\begin{thm}\cite[Corollary 3.11]{JLC}\label{Chabert}
Let $K/\Q$ be a cyclic extension of number fields. Then
$$
    |Po(K)|=\begin{cases}
     \frac{\prod_p e(p)}{2[K:\Q]}\;\; \text{if $K$ is real and } N_{K/\Q}(\mathcal{O}_K^\times) =\{1\}; \\ \\
     \frac{\prod_p e(p)}{[K:\Q]}\; \text{else}.
\end{cases}
$$
Here, the product runs over all rational primes $p$ ramifying in $K$ and $e(p)$ is the ramification index of $p$ in $K/\Q$.
\end{thm}

For cyclic number fields $K$ of odd degree, or for imaginary cyclic number fields $K$, from Theorem \ref{Chabert} we have
\begin{equation}\label{eC}
|Po(K)|=\frac{\prod_p e(p)}{[K:\Q]}.
\end{equation}
 In \cite{ZX85}, Zhang showed that for abelian number fields the narrow genus field can be easily constructed. More precisely for abelian number fields $K$ of degree $\ell^a$, for some prime $\ell$, the narrow genus field is described as given below
\begin{equation}\label{eZ1}
K^{*}=K \prod_{p \neq \ell} C_p = \prod_p C_p ~,
\end{equation}
where $p$ runs over rational primes which ramify in $K$, $e(p)$ is the ramification index of $p$ in $K$, $C_p$ is the unique subfield of degree $e(p)$ of $\Q(\zeta_p)$ whenever $p \neq \ell$ and $C_{\ell}$ is a subfield of degree $e(\ell)$ of $\Q(\zeta_{\ell^b})$ for some integer $b$.  Further, if $L$ is an abelian field which is compositum of $L_1$ and $L_2$ then $L^*=L_1^*L_2^*$. We recall that, as in Zhang \cite{ZX85}, the narrow genus field of an abelian field $K$ is maximal abelian extension $K^*/\Q$ containing $K$ such that $K^*/K$ is unramified at all finite primes. 

\begin{proof}
[Proof of Theorem $\ref{C1}$]
We assert that the extension $K^*/K$ is unramified at all primes (both finite and infinite).  Any cyclic field of degree $n$ is a compositum of cyclic fields of prime power degrees. We prove the assertion by assuming that $K$ is a cyclic field of prime power degree $\ell^a$. 
It is easy to see that if $K \not \subset \mathbb{R}$, then $\tilde{K}=K^*$. Now, suppose $K \subset \mathbb{R}$, then $\ell$ is an odd prime. Suppose $\sigma: K \longrightarrow \mathbb{R}$ is a field embedding and $\tilde{\sigma}: K^* \longrightarrow \C$ is an extension of $\sigma$.
First we note that in the description of $K^*$ as in (\ref{eZ1}), for each prime $p$ ramifying in $K$ we have $e(p)= \ell^c$, for the odd prime $\ell$ dividing $[K:\Q]$ and some integer $c>0$, and $C_p$ is the unique subfield of $\Q(\zeta_p)$ of degree $e(p)$. Note that $e(p)$ divides the degree $[\Q(\zeta_p)^{+}:\Q]$. Hence, we see that $C_p$ is a real Galois field. Consequently
$$\tilde{\sigma}\left(\prod_p C_p\right) \subset \mathbb{R}.$$
Also $\tilde{\sigma}(K) \subset \mathbb{R}$ and hence $\tilde{\sigma}(K^*) \subset \mathbb{R}$. This shows that the extension $K^*/K$ is unramified at infinite places as well. Therefore, we have $\tilde{K}=K^*$. Thus, $K^*/K$ is unramified at all primes in both cases: either when $K$ is a cyclic number field of odd degree, or when $K$ is an imaginary cyclic number field. \\
Now we suppose that $H(K)/\Q$ is abelian. Then $H(K)=K^*$, and from (\ref{eZ1}), we see that 
$$h_K=[H(K):K]=[K^*:K]=\frac{\prod_{p } e(p)}{[K:\Q]} .$$ 
From (\ref{eC}) it follows that $|Po(K)|=h_K$. Since $Po(K)$ is a subgroup of $C\ell(K)$, we conclude that $Po(K)=C\ell(K)$.\\
Conversely, we assume that $Po(K)=C\ell(K)$. From (\ref{eC}) and (\ref{eZ1}), we have $|Po(K)|=[K^*:K]$.  Thus $[K^*:K]=h_K$. Since $K^*/K$ is abelian and unramified at all primes, we have $K^* \subset H(K)$. We conclude that $H(K)=K^*$. Thus, $H(K)/\Q$ is abelian. 
\end{proof}

The following result is an immediate consequence of Theorem \ref{C1}.
\begin{cor}\label{C32}
  Let $K$  be an imaginary quadratic number field. Then, the extension $H(K)/\Q$ is abelian if and only if $Po(K)=C\ell(K)$.
\end{cor}

It is therefore natural to ask what happens in the case of real quadratic fields $\Q(\sqrt{d})$.  Our next results address this case.  Let $\epsilon_K$ denote the fundamental unit of the real quadratic field $K=\Q(\sqrt{d})$.

\begin{thm}\label{newrealquad}
Let $d > 0$ be a squarefree integer such that none of its prime divisors is congruent to $3 \pmod{4}$.  

\begin{enumerate}
  \item Assume that $N(\epsilon_K) \neq 1$. Then, $H(K)/\Q$ is an abelian extension if and only if $C\ell(K) \simeq {Po}(K)$.
  \item Assume that $N(\epsilon_K) = 1$. Then, $H(K)$ is an abelian extension of $\Q$ if and only if 
  $ [\,C\ell(K) : {Po}(K)\,] = 2$.
  
\end{enumerate}
\end{thm}
\begin{proof}
Let $d = p_1p_2\cdots p_r$, where the primes $p_i$ are distinct and satisfy 
$p_i \not\equiv 3 \pmod{4}$ for all $i$. Then, exactly $r$ many primes are ramified in $K$, and hence $[K^* : K] \;=\; 2^{r-1}$.
Now consider the multiquadratic field
\[
   L \;=\; \Q(\sqrt{p_1}, \sqrt{p_2}, \dots, \sqrt{p_r}).
\]
Clearly $L/K$ is unramified at all primes of $K$, and $L/\Q$ is an abelian extension of degree $2^r$. 

\noindent
We first prove ({\it{1}}). Since $K/\Q$ is a cyclic extension and $N(\epsilon_K) \neq 1$, using Theorem \ref{Chabert} we  have
\begin{equation}\label{realeq2}
   |Po(K)| \;=\; 2^{r-1}.
\end{equation}
Assume that $H(K)/\Q$ is abelian. Then 
\[
   H(K) \;\subseteq\; K^*.
\]
Moreover, $L \subseteq H(K) \subseteq K^*$ and since 
\[
   [L:\Q] = 2^r = [K^*:\Q],
\]
we must have $ L = H(K).$ Therefore,
\[
   h_K = [H(K):K] = [L:K] = 2^{r-1}.
\]
Combining this with \eqref{realeq2}, and noting that $Po(K)\subseteq C\ell(K)$, we obtain
\[
   C\ell(K) \;\simeq\; Po(K).
\]
Conversely, suppose $C\ell(K)\simeq Po(K)$. Then, by\eqref{realeq2}, 
$ h_K = 2^{r-1}.$
Since $[L:K]=2^{r-1}=h_K$, we conclude that $L=H(K)$. Therefore $H(K)/\Q$ is abelian. This proves {\it(1)} of Theorem~\ref{newrealquad}.\\

\noindent
Now we prove {\it (2)}.
Since $N(\epsilon_K)=1$,  from Theorem \ref{Chabert} we have
\[
   |Po(K)| = 2^{r-2}.
\]
Suppose that $H(K)/\Q$ is an abelian extension. Then, we again have
\[
   L \subseteq H(K) \subseteq K^*.
\]
Now,
\[
   [L:K] = 2^{r-1}, \qquad [K^*:K] = 2^{r-1}.
\]
Hence
\[
   2^{r-1} = [L:K] \leq [H(K):K] \leq [K^*:K] = 2^{r-1}.
\]
This forces
\[
   [H(K):K] = 2^{r-1}.
\]
Thus
\[
   h_K = |C\ell(K)| = [H(K):K] = 2^{r-1}.
\]
Since $|Po(K)|=2^{r-2}$ and $Po(K)\subseteq C\ell(K)$, it follows that
\[
   [\,C\ell(K):Po(K)\,] = \frac{|C\ell(K)|}{|Po(K)|} = \frac{2^{r-1}}{2^{r-2}} = 2.
\]
Conversely, suppose $[C\ell(K):Po(K)]=2$.  
Then
\[
   h_K = |C\ell(K)| = 2\cdot |Po(K)| = 2\cdot 2^{r-2} = 2^{r-1}.
\]
In particular, $h_K=[L:K]$, so we deduce $L=H(K)$.  
Since $L/\Q$ is abelian, it follows that $H(K)/\Q$ is abelian. This proves {\it(2)} of Theorem~\ref{newrealquad}.
\end{proof}
The case when there is a prime $p \equiv 3 \pmod{4}$ dividing $d$ is handled in the next theorem. 
\begin{thm}\label{realquad}
    Let $d > 0$ be such that there exists a prime $p \equiv 3 \pmod{4}$ dividing $d$. Then $H(K)/\mathbb{Q}$ is abelian if and only if ${Po}(K) = {C\ell}(K)$.
\end{thm}
\begin{proof}  Since $p \equiv 3 \pmod 4$, it follows that $\Q(\sqrt{-p}) \subset K^*$. On the other hand, $\Tilde{K} \subset \mathbb{R}$. Therefore, $K^* \neq \Tilde{K}$. Consequently, we have   $[K^* : \tilde{K}] = 2$. Let the discriminant $D_K$ have exactly $r$ distinct prime factors. Then,
\begin{equation*}\label{ke1}
    [K^* : K] = 2^{r-1}
\end{equation*}
and
\begin{equation}\label{ke2}
    [\tilde{K} : K] = 2^{r-2}.
\end{equation}
Now we claim that $N(\epsilon_K) = 1$. Suppose, for contradiction, that $N(\epsilon_K) = -1$, where $\epsilon_K = \alpha + \beta \sqrt{d}$ is the fundamental unit of $K$, with $\alpha, \beta \in \mathbb{Z}$. Then we have the negative Pell equation
\[
\alpha^2 - d\beta^2 = -1.
\]
Reducing both sides modulo $p$, we get $\alpha^2 \equiv -1 \pmod{p}$, which implies that $-1$ is a quadratic residue modulo $p$, that is, $\left(\frac{-1}{p}\right) = 1$. However, since $p \equiv 3 \pmod{4}$, we have
\[
\left(\frac{-1}{p}\right) = (-1)^{\frac{p-1}{2}} = -1,
\]
a contradiction. Hence, $N(\epsilon_K) = 1$. Now, applying Theorem \ref{Chabert}, we have
\begin{equation}\label{ke3}
    |Po(K)| = 2^{r-2}.
\end{equation}

Assume that $H(K)/\Q$ is an abelian extension. Then, $H(K) = \tilde{K}$. Therefore, using \eqref{ke2} and \eqref{ke3}, we conclude that $ C\ell(K) = Po(K)$.

Conversely, let  $ C\ell(K) = Po(K)$. Then, by \eqref{ke3} we have $h_K = 2^{r-2}$. Therefore, from \eqref{ke2} and the fact that $\tilde{K} \subset H(K)$, it follows that that $H(K) = \tilde{K}$, and hence $H(K)/\Q$ is an abelian extension. This completes the proof.
\end{proof}

The next result gives a very simple description of class groups of number fields of prime degree and with absolutely abelian Hilbert class field.
\begin{cor}
Let $K/\Q$ be a cyclic extension of degree $q$ for an odd prime $q$. If $H(K)/\Q$ is abelian then 
$$C\ell(K) \cong \left(\Z/q\Z\right)^{s-1},$$
where $s$ is the number of primes ramified in $K$.
\end{cor}
\begin{proof}
From the Corollary \ref{AC}, it follows that $q$ is the only prime factor of $|C\ell(K)|$. Consequently, we see that $$K^*=\widetilde{K}.$$
Now, the corollary follows from Theorem \ref{C1}, Theorem \ref{Chabert} and (\ref{eZ1}).
\end{proof}

In the remaining part of this section, we show that the assumption $H(K)/\Q$ is abelian imposes an upper bound on the class number of the Hilbert class field $H(K)$. For this, we need the following definition (see \cite{BR63}).
\begin{defn}
For any number field $E$ whose degree is divisible by $p^a$ for some prime $p$, we define 
$$R(E,p^a)=[E:\Q]\left(\frac{1}{p} + \ldots +\frac{1}{p^a} \right)+a.$$
For any two positive integers $x$ and $y$, we define
$$R(x,y)=\prod_p p^{R(E,p^{a(p)})},$$
where $y=\prod_p p^{a(p)}$ is the factorisation of $y$ and $E$ is any field of degree $xy$.
\end{defn}
We prove the following result.
\begin{thm}\label{N1}

    Let $K$ be a number field of degree $n$ with class number $h_K=h$. Let the Hilbert class field $H(K)$ over $\mathbb{Q}$ is abelian. Let $m\in\mathbb{Z}^+$ be the conductor of $H(K)$ and $[\mathbb{Q}(\zeta_m): H(K)]=m_1$. Then we have
    $$h_{H(K)}\leq \frac{R(nh,m_1)R(1,n)|Po(K)| |Po(\Q(\zeta_m)/H(K))|}{\varphi(m)}.$$
Here $Po(\Q(\zeta_m)/H(K))$ denotes the relative P\'olya group as considered by Chabert \cite{JLC1}.
\end{thm}

We now recall some results which are needed to prove Theorem \ref{N1}. The following result, describing the size of the first cohomology group, was obtained in \cite{JLC1}.
\begin{thm}\cite[Proposition 4.4]{JLC1}\label{thm1}
Let $K/F$ be a Galois extensions of number fields. Then
 $$|H^1(G_{K/F}, \mathcal{O}_K^\times)|=\frac{h_F \prod_{\mathfrak{p}}e(\mathfrak{p})}{|Po(K/F)|} ,$$
where the product runs over all the ramified prime ideals $\mathfrak{p}$ of $F$.
\end{thm}
The following result of Brumer and Rosen \cite{BR63} provides an upper bound on the size of the first cohomology group $H^1(G_{K/F},\mathcal{O}_K^\times)$ in terms of the parameters $R(L,p^t)$.
\begin{thm}\cite[Proposition 3.4]{BR63}\label{thm3}
Let $K/F$ be a Galois extension of number fields of degree $n$.  Suppose $n=\prod_p p^{a(p)}$ is the prime factorization of $n$ into prime powers. Then $|H^1(G_{K/F},\mathcal{O}_K^\times)|$ divides $\prod_p p^{R(L,p^{a(p)})}$.
\end{thm}
Next we mention the following result due to Zantema (see \cite[Page 9]{ZAN}) which connects the first cohomology group $H^1(G_K,\mathcal{O}_K^{\times})$ with the primes ramifying in $K/\mathbb{Q}$ and the P\'olya group $Po(K)$.
\begin{prop}\label{EXACT}
\cite[Theorem 1.3]{ZAN} Let $K/\mathbb{Q}$ be a Galois number field. Then there exists a canonical embedding $\displaystyle\Psi : H^1(G_K,\mathcal{O}_K^{\times})\longrightarrow \bigoplus_{p }\mathbb{Z}/e(p)\mathbb{Z}$ and an exact sequence of abelian groups
$$1\rightarrow H^1(G_K,\mathcal{O}_K^{\times})\xrightarrow{\Psi} \bigoplus_{p }\mathbb{Z}/e(p)\mathbb{Z}\rightarrow Po(K) \rightarrow1$$
where the direct sum runs over all rational primes ramifying in $K$.
\end{prop}
\begin{proof} [Proof of Theorem $\ref{N1}$]
We have $H(K)/\Q$ is abelian and $m$ is the conductor of $H(K)$. Thus we have the the following tower of field extensions.
\[
\begin{tikzcd}[row sep=1cm]
    & \mathbb{Q}(\zeta_m) \\
    & H(K) \ar[-, "m_1"]{u} \\
    & K \ar[-, "h"]{u} \\
    & \mathbb{Q} \ar[-, "n"]{u} \tikz[overlay,remember picture] \coordinate (bottom);
\end{tikzcd}
\]   
In this tower of fields, ramification happens only for the extensions $K/\Q$ and $\Q(\zeta_m)/H(K)$ and we have the following identity
\begin{equation}\label{eqN1}
\prod_p e(p) \prod_{\mathfrak{P}}e(\mathfrak{P})=\varphi(m),
\end{equation} 
where the first product runs over the rational primes $p$ ramifying in $K$ and the second product runs over the prime ideals $\mathfrak{P}$ of $H(K)$ ramifying in $\Q(\zeta_m)$.
 Using Theorem \ref{thm1} for the extension $\Q(\zeta_m)/H(K)$, we obtain
\begin{equation}\label{eqN2}
    \prod_{\mathfrak{P}}e(\mathfrak{P})=\frac{|H^1(\Gal{\mathbb{Q}(\zeta_m)/ H(K)}, \mathbb{Z}[\zeta_m]^{\times})| |Po(\Q(\zeta_m)/H(K))|}{h_{H(K)}},
\end{equation}
where the product runs over all prime ideals $\mathfrak{P}$ of $H(K)$ ramifying in the extension $\Q(\zeta_m)/H(K)$. Suppose $m_1=\prod_q q^{b(q)}$ is the prime factorization of $m_1$. Then, using the upper bound on $ |H^1(\Gal{\mathbb{Q}(\zeta_m)/ H(K)}, \mathbb{Z}[\zeta_m]^{\times})|$ from Theorem \ref{thm3}, we see that
\begin{equation}\label{eqN3}
 \prod_{\mathfrak{P}}e(\mathfrak{P})  \leq \frac{\prod_q q ^{R(\mathbb{Q}(\zeta_m),q^{b(q)})} |Po(\Q(\zeta_m)/H(K))| }{h_{H(K)}}=\frac{R(nh,m_1) |Po(\Q(\zeta_m)/H(K))| }{h_{H(K)}}.
\end{equation}
On the other hand, from Proposition \ref{EXACT} we get
\begin{equation}\label{eqN4}
  \prod_p e(p)= |H^1(G_K,\mathcal{O}_K^\times))|\cdot |Po(K)|,
\end{equation}
where the product runs over all the rational primes ramifying in $K$. Again, using the bound on $  |H^1(G_K,\mathcal{O}_K^\times))|$ from Theorem \ref{thm3} we get
\begin{equation}\label{eqN5}
     \prod_pe(p)\leq R(1,n)\cdot |Po(K)|.
\end{equation}
From equations (\ref{eqN1}), (\ref{eqN3}) and (\ref{eqN5}) it follows that
\begin{equation}\label{eq6}
  h_{H(K)}\leq \frac{R(nh,m_1)R(1,n) |Po(K)| |Po(\Q(\zeta_m)/H(K))|}{\varphi(m)}.
\end{equation}
\end{proof}

\section{Proof of Theorem \ref{Main}}\label{secC}
In this section, we prove a stronger theorem than Theorem \ref{Main}. Before stating the theorem we recall that Hilbert $\ell-$class field $H_{\ell}(K)$ of $K$ is the largest subfield $H_{\ell}(K)$ of $H(K)$ such that the degree $[H_{\ell}(K):K]$ is a power of $\ell$.
\begin{thm}\label{Main1}
Let $K$ be an abelian number field of degree $n$ and $\ell$ be a prime number. Suppose that the Hilbert $\ell-$class field $H_{\ell}(K)$ of $K$ is an abelian extension of $\Q$. Then for any positive $\epsilon$, we have
\begin{equation}\label{ME2}
|C\ell(K)[\ell]| \ll_{\epsilon,\ell, n} D_K^{\epsilon},
\end{equation}
where $D_K$ is the absolute value of the discriminant and the constant depends only on $n, \ell$ and $\epsilon$.
\end{thm}
To prove Theorem \ref{Main1}, we begin with an analog of Corollary \ref{AC}. This can be proved analogous to Theorem \ref{OT}, but for the sake of completeness we include a proof.
\begin{prop}\label{AC1}
Let $K$  be an abelian number field of degree $n$. If there exists a prime $\ell$ such that $\ell$ divides $h_K$ but does not divide $n$, then the extension $H_{\ell}(K)/\Q$ is not abelian.
\end{prop}
\begin{proof}

Let $\sigma$ be a non-trivial automorphism in $G_K$, and $\widetilde{\sigma} \in G_{H_{\ell}(K)}$ be an extension of $\sigma$. Let $\mathfrak{A}$ be an element of $C\ell(K)$ of order $\ell$. Theorem \ref{CCT} implies that $\mathfrak{A}$ contains a prime ideal $\mathfrak{p}$ of $K$ of residue degree one, such that $\mathfrak{p}$ is unramified in $K/\Q$. Let 
\[
\tau := \left( \frac{\mathfrak{p}}{H_{\ell}(K)/K} \right)
\]
 be the Frobenius automorphism of $\mathfrak{p}$. By Theorem \ref{HCF}, the order of $\tau$ is $\ell$. Since $\tau \in G_{H_{\ell}(K)/K}$, it follows that $\tau$ is an element of $G_{H_{\ell}(K)}$.  For brevity, denote $\text{Frob}_\mathfrak{p} := \text{Frob}_\mathfrak{P}=\tau$, where $\mathfrak{P}$ is a prime of $H_{\ell}(K)$ above $\mathfrak{p}$. By definition, for any $x \in \mathcal{O}_{H_{\ell}(K)}$,
\[
\text{Frob}_\mathfrak{P}(x) \equiv x^{N(\mathfrak{p})} \pmod{\mathfrak{P}},
\]
where $N(\mathfrak{p}) = |\mathcal{O}_K/\mathfrak{p}|$. For $x\in\mathcal{O}_{H_{\ell}(K)}$ we have $\widetilde{\sigma}^{-1}(x)\in\mathcal{O}_{H_{\ell}(K)}$. Therefore
\begin{equation*}
    \text{Frob}_\mathfrak{P}(\widetilde{\sigma}^{-1}(x))\equiv (\widetilde{\sigma}^{-1}(x))^{N(\mathfrak{p})} \pmod{\mathfrak{P}}.
\end{equation*}
Applying $\widetilde{\sigma}$ to this congruence yields
\begin{equation*}
    \left(\widetilde{\sigma} \circ \text{Frob}_\mathfrak{P} \circ \widetilde{\sigma}^{-1}\right)(x) \equiv x^{N(\mathfrak{p})} \pmod{\widetilde{\sigma}(\mathfrak{P})}.
\end{equation*}
Thus,
\begin{equation*}
    \text{Frob}_{\mathfrak{\widetilde{\sigma}(P)}} (x) \equiv x^{N(\mathfrak{p})} \pmod{\widetilde{\sigma}(\mathfrak{P})}.
\end{equation*}
Since $N(\mathfrak{p})=N(\sigma(\mathfrak{p}))$, we have
\begin{equation*}\label{N11}
    \text{Frob}_{\mathfrak{\widetilde{\sigma}(P)}} (x) \equiv x^{N(\sigma(\mathfrak{p}))} \pmod{\widetilde{\sigma}(\mathfrak{P})}.
\end{equation*}
Since $\widetilde{\sigma}$ is a lift of $\sigma$, it maps primes of $H_{\ell}(K)$ above $\mathfrak{p}$ to those above $\sigma(\mathfrak{p})$. In particular, $\widetilde{\sigma}(\mathfrak{P})$ lies above $\sigma(\mathfrak{p})$.  Thus,
\begin{equation}\label{Nimp}
\widetilde{\sigma} \circ \text{Frob}_\mathfrak{p} \circ \widetilde{\sigma}^{-1} = \text{Frob}_{\sigma(\mathfrak{p})}.    
\end{equation}
As $H_{\ell}(K)/\Q$ is abelian, using (\ref{Nimp}) we get
$$\left( \frac{\sigma(\mathfrak{p})}{H_{\ell}(K)/K} \right)=\tau ,\mbox{ for all }\sigma \in G_K.$$
Consequently,
\begin{equation*}
    \displaystyle\prod_{\sigma\in G_K}\left( \frac{\sigma(\mathfrak{p})}{H_{\ell}(K)/K} \right) =\tau^n.
\end{equation*}
On the other hand, from the multiplicativity of Frobenius elements, we have
\begin{equation*}\label{Neq1}
    \displaystyle\prod_{\sigma\in G_K}\left( \frac{\sigma(\mathfrak{p})}{H_{\ell}(K)/K} \right) = \left( \frac{N(\mathfrak{p})}{H_{\ell}(K)/K} \right)=Id,
\end{equation*}
as $N(\mathfrak{p}) = \prod_{\sigma \in G_K}\sigma(\mathfrak{p})$ is a principal ideal. 
Thus, $$\tau^n=Id.$$
This proves that $\ell |n$.
\end{proof}
\begin{proof}[Proof of Theorem \ref{Main1}]
If $C\ell(K)[\ell]$ is trivial then theorem holds trivially. We assume that $C\ell(K)[\ell]$ is non-trivial. From Proposition \ref{AC1} it follows that $\ell|[K:\Q]$. Let $[K:\Q]=\prod_{i=1}^t \ell_i^{e_i}$ with $\ell_1=\ell$ and $e_1=e$.  As $K/\Q$ is abelian, for each $i$, there is unique subfield $K_i$ of $K$ of degree $\ell_i^{e_i}$. We have $K_1 \ldots K_t=K$ and $K_i \cap K_j=\Q$ for all $i \neq j$. We put 
$$F=K_2 \ldots K_t.$$ 
Then $K_1 \cap F=\Q$.\\

Suppose $G(K)$ is the genus field of $K$. As $H_{\ell}(K)/\Q$ is abelian, we see that $H_{\ell}(K) \subset G(K)$. Consequently,
$$[H_{\ell}(K):K]|g_K, \mbox{ where }g_K=[G(K):K] \mbox{ is the genus number of }K.$$
It is well known that (for example see page 56 of Ishida's book)
$$g_K=\frac{\prod e_p(K/\Q)}{[K:\Q]}, $$
where the product runs over all the primes $p$ ramifying in $K$ and $e_p(K/\Q)$ is the ramification index of $p$ in $K/\Q$. Let $e_p(F/\Q)$ and $e_p(K/F)$ denote the ramification index of $p$ for the extension $F/\Q$ and $K/F$ respectively. Then
$$g_K= \frac{\prod e_p(K/F) \prod e_p(F/\Q)}{[K_1:\Q][F:\Q]}.$$
Note that $\ell \nmid [F:\Q]$, hence we get
$$[H_{\ell}(K):K]| \frac{ \prod e_p(K/F)}{[K_1:\Q]}.$$
Furthermore, $e_p(K/F)=e_p(K_1F/F)=e_p(K_1/\Q)$. Thus 
$$[H_{\ell}(K):K] | \frac{\prod e_p(K_1/\Q)}{[K_1:\Q]}.$$
Suppose $D_K$ has $s$ distinct prime factors then
$$\frac{\prod e_p(K_1/\Q)}{[K_1:\Q]} \leq [K_1:\Q]^{(s-1)}=\ell^{e(s-1)}.$$
Combining these we get 
$$[H_{\ell}(K):K] \leq \ell^{e(s-1)},$$
and $|C\ell(K)[\ell]| \leq \ell^{e(s-1)}.$ From this, it is immediate to see that
$$|C\ell(K)[\ell]| \ll_{\epsilon, \ell,n}D_K^{\epsilon}.$$
\end{proof}

Now we intend to show that there are infinitely many number fields $K$ whose Hilbert $\ell-$class field $H_\ell(K)$ is an abelian extension of $\Q$. We assume that $\ell$ is an odd prime number. For any non-negative integer $s$ and positive real number $x$, we let $M_{s,x}$ denote the number of cyclic extensions $K/\Q$ of degree $\ell$ satisfying the following conditions:\\
i) the conductor $f_K$ of $K$ is divisible by exactly $s$ distinct prime numbers; ii) $f_K \leq x$; iii) the $\ell-$part of class number of $K$ satisfies $h_{\ell}(K)=\ell^{s-1}$. Furthermore, let $N_{s,x}$ denote the number of cyclic extensions $K/\Q$ of degree $\ell$ with $h_{\ell}(K)=\ell^s$ and with $f_K \leq x$. In \cite{FG83} (see Theorem 1), the author established the following asymptotic.

\begin{equation}\label{FrankG}
N_{s,x}=M_{s+1,x}+ O\left(\frac{M_{s+1,x}}{\log \log x}\right) \mbox{ as }x \longrightarrow \infty.
\end{equation}
From (\ref{FrankG}), it follows that a positive proportion of cyclic number fields $K$ of degree $\ell$ over $\Q$ with $h_{\ell}(K)=\ell^s$ have exactly $s+1$ ramified primes. Now, from (\ref{eZ1}) it follows that $[K^*:K]=\ell^s$. Consequently, we see that $K^*=H_{\ell}(K)$ for a positive proportion of cyclic number fields $K$ of degree $\ell$ over $\Q$ with $h_{\ell}(K)=\ell^s$. This establishes the following.
\begin{thm}\label{infinite}
For a positive proportion of cyclic extensions $K/\Q$ of degree $\ell$, the Hilbert $\ell-$class field $H_{\ell}(K)$ is an abelian extension of $\Q$.
\end{thm}
As a consequence of Theorem \ref{infinite}, we know that there are infinitely many fields for which Theorem \ref{Main1} and Theorem \ref{Main} are true. As mentioned in the introduction, the $\ell-$torsion conjecture was recently proved for $\ell-$extensions of $\Q$. The additional feature of Theorem \ref{Main1} is that it holds for all extensions $K/\Q$ for which $H_{\ell}(K)$ is an abelian extension of $\Q$, irrespective of whether $K/\Q$ is an $\ell-$extension or not. It will be interesting to see if there are abelian extensions $K/\Q$ which are not $\ell-$extensions but  $H_{\ell}/\Q$ is abelian and $C\ell(K)[\ell]$ is non-trivial. Also, in light of Theorem \ref{Main1}, Theorem \ref{infinite}, and the works of Kl\"{u}ners and Wang \cite{KW22} it looks like that a complete proof of $\ell-$torsion conjecture will need a much deeper insight.\\

\section{Some finiteness results and concluding remarks}\label{secE}

In this section we show that there are very few imaginary abelian number fields with absolutely abelian Hilbert class field. We recall the following result.
\begin{thm} \cite[Theorem 1]{MR1416716} \label{NT1}
 Let \(\epsilon > 0\) be given.  There exists \(d_\epsilon\) such that for any imaginary abelian number field \(K\) satisfying $d_{F} > d_{\epsilon}$ we have $g_{F} < \bigl(h_{F}^-\bigr)^{\epsilon}, $ where $d_F$, $g_F$  and $ h_K^{*}$ denote the absolute value of the discriminant of $K$, the narrow genus number of $K$ and the relative class number of $K$, respectively.
Therefore, there are only finitely many imaginary abelian number fields \(K\) such that \(h_{K} = g_{}\).  Moreover, there exists \(d_{1}\), effectively computable, such that $h_{K} = g_{K}$ implies $d_{K} \le d_{1}$, 
provided that \(K\) is neither an imaginary quadratic number field nor an imaginary biquadratic bicyclic number field.
\end{thm}

For imaginary number fields, the narrow genus field and the absolute genus field are same. The following result can be deduced immediately from Theorem~\ref{NT1}.

\begin{thm}
Only finitely many imaginary abelian number fields have absolutely abelian Hilbert class fields.
\end{thm}

Now, we recall the following results describing imaginary abelian number of fields for which genus number is same as the class number.

\begin{thm} \cite[Theorem 2]{MR1362937} \label{NT2} If the Generalized Riemann Hypothesis is true, we have exactly $301$ imaginary abelian number fields $K$ of type $(2,2,\dots,2)$, i.e. $K=\mathbb{Q}(\sqrt{-m_1}, \sqrt{-m_2},\dots, \sqrt{-m_r})$, such that $[C\ell(K): G_{\widetilde{K}/K}]=1$. Among them 
\begin{enumerate}
    \item $65$ fields are quadratic;
    \item $219$ fields are bicyclic biquadratic;
    \item $17$ fields are of type $(2,2,2)$ (Explicit fields are listed in the article \cite{MR1362937}).
\end{enumerate}
\end{thm}
\begin{thm}\cite{MR1823189}\label{NT3}
There are exactly $424$ imaginary non-quadratic abelian number fields with class numbers equal to their genus class numbers. Among them:
\begin{enumerate}
    \item $77$ are cyclic,
    \item $347$ are non-cyclic.
\end{enumerate}
Their degrees are less than or equal to $24$, their class numbers are $1$, $2$, or $4$. The conductors of these fields are less than or equal to $65689$.

\end{thm}

Using \thmref{NT2} and \thmref{NT3}, we can also deduce the following result, which offers a complete classification of imaginary abelian number fields with absolutely abelian Hilbert class fields.
\begin{thm}
Let $K$ be an imaginary abelian number field. Then:
\begin{enumerate}
    \item If $K$ is quadratic, i.e., $K = \mathbb{Q}(\sqrt{-d})$, then, under the generalized Riemann Hypothesis, there are exactly $65$ values of $d$ for which the extension $H(K)/\mathbb{Q}$ is abelian.

    \item If $K$ is not quadratic, then there are exactly $424$ such fields for which the extension $H(K)/\mathbb{Q}$ is abelian. This result is unconditional.
\end{enumerate}
\end{thm}

From the Cohen--Lenstra heuristics, we expect that about $75\%$ of real quadratic fields with prime discriminant have class number 1, and hence have an absolutely abelian Hilbert class field.  However, no infinite family of number fields with absolutely abelian Hilbert class field is known. Below we provide families of real quadratic and real quartic fields whose Hilbert class field is not absolutely abelian.
\begin{thm}\label{EXT2}
 Let $p=b^2+c^2\equiv 1 \pmod4$ be a prime, where $b$ and $c$ are positive integers. Suppose $K_{a,p}=\Q(\sqrt{a(p+b\sqrt{p})})$, where $a$ is a square-free odd integer such that $|a|>1$ and is relatively prime to $p$. If the class number of $\Q(\sqrt{p})$ is greater than one, then the Hilbert class field of $K_{a,p}$ can not be absolutely abelian. Moreover, the Hilbert class field of $\Q(\sqrt{p})$ can not be absolutely abelian.  
\end{thm}
\begin{proof}
Since $p \equiv 1 \pmod{4}$,  from genus theory, it is known that the 2-rank of the class group of $\mathbb{Q}(\sqrt{p})$ is 0. However, we have assumed that $h_{\mathbb{Q}(\sqrt{p})} > 1$. This implies that there exists an odd prime $q$ such that $q \mid h_{\mathbb{Q}(\sqrt{p})}$. Now, by using Corollary \ref{AC}, we conclude that the Hilbert class field of $\mathbb{Q}(\sqrt{p})$ cannot be absolutely abelian.

Let us consider the cyclic quartic field $K_{a,p}=\mathbb{Q}(\sqrt{a(p+b\sqrt{p})})$. Therefore, we have $[\mathbb{Q}(\sqrt{a(p+b\sqrt{p})}) :\Q(\sqrt{p})]=2$. Thus, any odd prime that divides $h_{\mathbb{Q}(\sqrt{p})}$ also divides $h_{K_{a,p}}$, in particular $q\mid h_{K_{a,p}}$. Again using Corollary \ref{AC}, we conclude that the Hilbert class field of $K_{a,p}$ cannot be absolutely abelian.
\end{proof}

Now, we give a real quadratic field $K$ with absolutely abelian Hilbert class field $H(K)$ for which $Po(K) \neq C\ell (K)$. For this, let $K=\Q(\sqrt{221})$. We find that $h_K=2$, and hence $H(K)$ is absolutely abelian. However, as shown in \cite{ZAN}, $K$ is a P\'olya field. Thus $Po(K) \neq C\ell(K)$.

In the remaining part of this section we present our investigation of real quadratic fields $\mathbb{Q}(\sqrt d)$ with discriminant $\leq 10^5$ and whose Hilbert class field is absolutely abelian.  By Theorem~\ref{QN}, these are precisely the real quadratic fields whose class group is trivial or an elementary abelian $2$-group. Table~\ref{tab:counts} gives the distribution of such fields according to the $2$-rank of the class group.

\begin{table}[htbp]
\centering
\caption{Real quadratic fields of discriminant $\leq10^5$ with class group $(\mathbb{Z}/2\mathbb{Z})^k$.}
\label{tab:counts}
\begin{tabular}{c c l}
\toprule
$k$ & Count $\mathbb{Q}(\sqrt d)$ & Examples of $d$ \\
\midrule
0 & 10\,079 & $5,\ 8,\ 12,\ 13,\ 17,\ 21,\ 24,\ 28,\ \dots$ \\
1 &  8\,916 & $40,\ 60,\ 65,\ 85,\ 104,\ 105,\ 120,\ \dots$ \\
2 &  2\,901 & $520,\ 680,\ 780,\ 840,\ 924,\ 1020,\ \dots$ \\
3 &    287 & $4620,\ 7980,\ 8840,\ 9240,\ \dots$ \\
4 &      2 & $60060,\ 78540$ \\
$\ge 5$ & 0 & \\
\bottomrule
\end{tabular}
\end{table}

The data were obtained using the LMFDB database of number fields \cite{LMFDB}, restricting to real quadratic fields of discriminant $\leq 10^5$. All fields were filtered according to their class group structure. The complete list is omitted for brevity.

\section{Primes of higher degree} \label{secD}

In connection with Question 2, we consider the following set for any extension $K/F$ of number fields.

\begin{align*}
\mathbf{R}_{K/F}=\left\{\begin{array}{l}
f\in\mathbb{N}: \text{ $C\ell(K)$}\; \;\text{is generated by the ideal classes containing  }\\ \\ \qquad \qquad\text {  ideals  of}
     \; K \;\text{ unramified prime of residue degree}\; f \end{array}  \right\}.
\end{align*} 

It is quickly seen from Theorem \ref{CCT} that $1\in \mathbf{R}_{K/F}$. Before stating the main result of this section we recall the following group-theoretical lemma (for example, see \cite{PM22}).
\begin{lem}\label{HL28}
Let $G_1$ and $ G_2$ be two finite groups and let $\psi: G_2 \longrightarrow Aut(G_1)$ be a group homomorphism. Assume that $g_1 \in G_1$ and $ g_2 \in G_2$ have orders $n_1$ and $n_2$, respectively. If $\psi(g_2)$ is identity automorphism, then the order of the element $(g_1, g_2) $ in the semidirect product $G_1 \rtimes_{\psi} G_2$ is $lcm~ (n_1, n_2)$.
\end{lem}

The main result we prove in this section is the following theorem.
\begin{thm}\label{PHT1}
 Let $K/\Q$ be a Galois extension of degree $n$ and $h_K=u$. Assume that $f>1$ is a fixed divisor of $n$ such that $(f, u)=1$. Suppose $\psi: G_K \rightarrow \text{Aut}(G_{H(K)/K})$ is a group homomorphism such that
 $$ G_{H(K)} \simeq G_{H(K)/K} \rtimes_\psi G_K, $$ and $G_K$ has an element $\sigma$ of order $f$ such that $\psi(\sigma)=Id$. If $G_K$ has no element of order greater than $f$, then $f\in \R_{K/\Q}$.
 \end{thm}   



\begin{proof}
As $C\ell(K)$ is an abelian group and $C\ell(K) \cong G_{H(K)/K}$, it follows that there exist positive divisors $n_1, \dots, n_t$ of $u$ such that 
\begin{equation}\label{hpE1}
  G_{H(K)/K} \cong \Z/n_1\Z \oplus \dots \oplus \Z/n_t\Z. 
\end{equation} From \eqref{hpE1}, we find that there exist elements $\sigma_1, \dots, \sigma_t  $ in $G_{H(K)/K}$ of orders $n_1, n_2, \ldots, n_t$ respectively  such that 
\begin{equation}\label{hpE2}
G_{H(K)/K}=\left\{\prod_{i=1}^t \sigma_i^{j_i}: 1 \leq j_i \leq n_i , \forall i\right\}.
\end{equation}
From \eqref{hpE2} and the assumption that $f$ is relatively prime to $u$, we conclude that $G_{H(K)/K}= G_{H(K)/K}^f$, i.e.,
\begin{equation}\label{hpE3}
G_{H(K)/K}=\left\{\prod_{i=1}^t \sigma_i^{fj_i}: 1 \leq j_i \leq n_i , \forall i\right\}.
\end{equation}
By hypothesis, there exists an element $\sigma \in G_K$ of order $f$ such that $\psi(\sigma)=Id$. From Lemma \ref{HL28}, the order of the element $(\sigma_i, \sigma)$ is $n_i\cdot f$ for each $i$. As $\Psi$ is an isomorphism, there are elements $\tau_i \in G_{H(K)}$ such that $\Psi(\tau_i)=(\sigma_i, \sigma)$. Applying  Theorem \ref{CDT}, we obtain unramified prime ideals $\mathfrak{P}_i$ of $\mathcal{O}_{H(K)}$ such that
\begin{equation*}
\left(\frac{\mathfrak{P}_i}{H(K)/\Q} \right)= \tau_i \mbox{ for each }i=1, \dots, t.
\end{equation*} 
Let $\mathfrak{p}_i$ and $p_i$ denote the primes below $\mathfrak{P}_i$ of the fields $K$ and $\Q$, respectively. We have the following identities:
$$f(\mathfrak{P}_i|p_i)= \text{order of}\;\tau_i= n_i\cdot f \mbox{ and } f(\mathfrak{P}_i|p_i)=f(\mathfrak{P}_i|\mathfrak{p}_i)\cdot f(\mathfrak{p}_i|p_i).$$
Note that $f(\mathfrak{P}_i|\mathfrak{p}_i)$ is a divisor of $[H(K):K]$, $f(\mathfrak{p}_i|p_i)$ is a divisor of $n$, $f$ is relatively prime to $[H(K) : K]$ and $G_K$ has no element of order greater than $f$. Thus, we conclude the following:
\begin{equation*}
f(\mathfrak{P}_i|\mathfrak{p}_i)= n_i \mbox{ and } f(\mathfrak{p}_i/p_i)=f.
\end{equation*}
Let $$\left(\frac{\mathfrak{P}_i}{H(K)/K} \right)= \tau_i^{'} \mbox{ for each }i=1, \dots, t.$$ 
Then, $\tau_i' \in G_{H(K)/K}$  and the $\text{order of}\;\tau_i'=n_i$ for each $i$. We claim that the Galois group $G_{H(K)/K}$ is generated by $\tau_1', \dots, \tau_t'$.  From Lemma \ref{FR}, we get $\tau_i^f=\tau_i^{'}$ for each $i$. Therefore,
$$\Psi(\tau_i'^{j_i})=\Psi(\tau_i^{fj_i})=(\sigma_i, \sigma)^{fj_i}.$$ Now, using Lemma \ref{HL28}, we have
\begin{equation}\label{hpME5}
\Psi(\tau_i'^{j_i})=(\sigma_i^{fj_i}, 1).
\end{equation} 
Since $G_{H(K)/K}$ embeds inside $G_{H(K)/K} \rtimes_{\psi} G_K$ via the map $g \mapsto (g,1)$ and $\Psi$ is an isomorphism, our claim follows from \eqref{hpE3} and (\ref{hpME5}).\\
As the extension $H(K)/K$ is abelian, we have
$$\left(\frac{\mathfrak{p}_i}{H(K)/K} \right)=\left(\frac{\mathfrak{P}_i}{H(K)/K} \right)=\tau_i'.$$
Using Theorem \ref{HCF}, we conclude that $C\ell(K)$ is generated by the ideal classes of the prime ideals $\mathfrak{p}_i$. Consequently, $f \in \mathbf{R}_{K/\Q}$, which completes the proof.
\end{proof}

It is well known that the vanishing of the second cohomology group $H^2(G_K,$ $ C\ell(K))$ implies
$G_{H(K)}\cong  G_{H(K)/K} \rtimes_{\psi} G_K$ (see \cite[\S 17.4]{DF04}). From the proof of Theorem \ref{PHT1}, it seems that number fields $K$ satisfying $H^2(G_K, C\ell(K))=0$ might be immediate target for examples of fields with non-trivial $\mathbf{R}_{K/\mathbb{Q}}$.  However, it is also possible to obtain number fields $K$ with non-trivial $\mathbf{R}_{K/\mathbb{Q}}$ even if $H^2(G_K, C\ell(K)) \neq 0$.  Our next result is in this direction.
\begin{thm}\label{PHT2} Let $u \geq 1$ be a square-free integer such that each prime divisor of $u$ is congruent to $2$ modulo $3$. Further, let $K/\Q$ be a Galois extension such that $G_K\simeq S_3$, where $S_3$ is the permutation group of order $6$. If $C\ell(K)\simeq \Z/{u\Z}$, then $3 \in\R_{K/\Q}$.
\end{thm}

\begin{proof}
Given that $C\ell(K) \simeq \Z/{u\Z}$. From Theorem \ref{HCF}, we have $G_{H(K)/K} \simeq C\ell(K)$. Consequently,  $G_{H(K)/K}\simeq \Z/{u\Z}$. Let $\ell$ be a prime divisor of $u$, and let $L$ be a subfield of $H(K)$ containing $K$ such that $[L : K]=\frac{u}{\ell}$. Using the facts that $H(K)/\Q$ is a Galois extension and $H(K)/K$ is an abelian extension, it is easy to see that $L$ is unique and $L/\Q$ is a Galois extension.  As $ L/\Q$ is a Galois extension, it follows that $G_{H(K)/L}$ is a normal subgroup of $G_{H(K)}$. Using the fundamental theorem of Galois theory, we obtain that the order of the $G_{H(K)/L}$ is $\ell$.

Since $3$ divides $[H(K) : \Q]$, it follows that $G_{H(K)/\Q}$ has a subgroup of order $3$, say $H$. Let $G=G_{H(K)/L}H$. We have 
\begin{equation}\label{hpE5}
    |G|= \frac{|G_{H(K)/L}| \cdot |H|}{|G_{H(K)/L}\cap H|}= \frac{3\ell}{|G_{H(K)/L}\cap H|}.
\end{equation}
Using the hypothesis that each prime divisor of $u$ is congruent to $2$ modulo $3$, we obtain that $|G_{H(K)/L} \cap H|=1$. Therefore, \eqref{hpE5} gives $|G|=3\ell$. Using the fact that $G_{H(K)/L}$ is a normal subgroup of $G_{H(K)}$, we conclude that $G$ is a subgroup of $G_{H(K)}$. Since $\ell$ is congruent to  $2$ modulo $3$, it implies that $H$ is a normal subgroup of $G$. Thus, we have 
$$G \cong G_{H(K)/L} \oplus H.$$ This concludes that $G_{H(K)}$ has an element $\tau$ of order $3\ell$. Theorem \ref{CDT} implies that there exists an unramified prime ideal $\mathfrak{P}$ in $H(K)$ with residue degree $3 \ell$ such that  
\begin{equation*}
\left(\frac{\mathfrak{P}}{H(K)/\Q} \right)= \tau.
\end{equation*} 
Let $\mathfrak{p}$ and $p$ denote the primes below $\mathfrak{P}$ of the fields $K$ and $\Q$, respectively. Note the following identities:
\begin{equation}\label{phe6}
    f(\mathfrak{P}|p)= \mathrm{ord}(\tau)= 3\ell \mbox{ and } f(\mathfrak{P}|p)=f(\mathfrak{P}|\mathfrak{p})\cdot f(\mathfrak{p}|p).
\end{equation}
We know that 
\begin{equation}\label{phe7}
    f(\mathfrak{P}|\mathfrak{p})\mid [H(K) : K]\; \; \text{and}\;\; f(\mathfrak{p}|p) \mid [K :\Q].
\end{equation}\\
{Case- (i):} Suppose $\ell$ is an odd prime. As $3$ does not divide $u$ and the degree of $K$ over $\Q$ is $6$, it follows that $\ell$ is relatively prime to $[K:\Q]$. Thus, using \eqref{phe6} and \eqref{phe7} we obtain 
 $$f(\mathfrak{P}|\mathfrak{p})=\ell\; \; \text{and}\;f(\mathfrak{p}|p)=3.$$ 
{Case- (ii):}  Suppose $\ell$ is an even prime. As $3$ is relatively prime to $u$ and $G_K$ has no element of order $6$, using \eqref{phe6} and \eqref{phe7}  we find that 
 $$f(\mathfrak{P}|\mathfrak{p})= \ell\; \; \text{and}\;f(\mathfrak{p}|p)=3.$$
 In both cases, we obtain $f(\mathfrak{P}|\mathfrak{p})=\ell\; \; \text{and}\;f(\mathfrak{p}|p) = 3.$
 Since $H(K)/K$ is an abelian extension, it follows that 
 $$\left(\frac{\mathfrak{P}}{H(K)/K} \right)= \left(\frac{\mathfrak{p}}{H(K)/K} \right)$$ and $\left(\frac{\mathfrak{p}}{H(K)/K} \right) \in G_{H(K)/K}$. Applying Theorem \ref{HCF}, we obtain that $C\ell(K)[\ell-]$ is generated by the ideal class of the prime ideal $\mathfrak{p}$. Consequently, $3  \in \mathbf{R}_{K/\Q}[\ell]$. This holds for each prime divisor $\ell$ of $m$. This concludes the proof.
\end{proof}

In the remaining part of this section, we aim to provide some examples using Theorem \ref{PHT2}. For this, we recall the following result from \cite{DF04}.
\begin{prop}\label{PHP26}\cite[Corollary $38$, Ch. $17$]{DF04} Let $A$ be  an abelian group and $G$ be a finite group such that $A$ is a $G$-module. If $|A|$  and $|G|$ are relatively prime, then every extension of $G$ by $A$ splits. 
\end{prop}

 Let $f(x) \in \Q[x]$ be a polynomial of degree $n$ with roots $\alpha_1, \alpha_2, \ldots, \alpha_n$. Then, the discriminant of $f$ is defined as $$D_f = \prod_{i<j}(\alpha_i -\alpha_j)^2.$$ 
The following lemma, from \cite{DF04},  provides the explicit formula for the discriminant of a cubic polynomial.
\begin{lem}\label{L41}\cite[Ch. $14$]{DF04}
Let $f(x)=x^3+ax^2+bx+c$ be a cubic polynomial. Then its discriminant $D_f$ in terms of $a, b$ and $c$ is given by $$D_f=a^2b^2-4b^3-4a^3c-27c^2+18abc.$$   
\end{lem}
Next, we provide some explicit families with Galois group $S_3$.
\begin{prop}\label{P42}
Let $f(x) = x^3+cx+c$  be a polynomial over $\mathbb{Z}$, where $c $ is an odd integer such that $c \equiv 1 \pmod 3$. Suppose $K$ is the splitting field of $f(x)$. Then, $G_{K} \simeq S_3$.
\end{prop}
\begin{proof}
Reading $f(x)$  modulo $2$, we obtain that $f(x)$ has no zeros in $\Z/{2\Z}$. Consequently, $f(x)$  is irreducible over $\Z$. Therefore, $f(x)$  is irreducible over $\Q$. Using Lemma \ref{L41}, we find that
\begin{align*}
    D_f &= -4c^3-27c^2 \\
   & = -c^2(4c+27).
\end{align*}
By hypothesis, $c \equiv 1 \pmod 3$, we conclude that $\sqrt{D_f}\not \in \Q$. Thus, $G_{K} \simeq S_3$. This completes the proof.
\end{proof}
Let $c=121$ in Proposition \ref{P42}. Using SageMath, we obtain that $h_K = 5$. Here, $[K :\Q]$ and $h_K$ are relatively prime. Therefore, the sequence 
\begin{equation*}
1 \longrightarrow G_{H(K)/K} \longrightarrow G_{H(K)} \longrightarrow G_{K} \longrightarrow 1
\end{equation*}
splits. As $3$ does not divide $|\text{Aut}(C\ell(K))|$, it follows that if $\psi : G_{K}\rightarrow \text{Aut}(C\ell(K))$ is any group homomorphism, there exits an element $\sigma \in G_{K}$ of order $3$ satisfying $\psi(\sigma) =Id$. We see that $G_{K}$ has no element of order greater than $3$. Thus, Theorem \ref{PHT1} applies and gives $3 \in \textbf{R}_{K/\Q}$.\\

For $c =  -59, -29, 11, 23, 25, 59, 71, 83$, we see that either $c \equiv 1 \pmod 3$, some values of $c$ or $\sqrt{D_f}\not \in \Q$. Therefore, $G_K \simeq S_3$. Using SageMath, we obtain that $h_K =2, 2, 2, 5, 2, 14, 14, 14 $ respectively. Here, $h_K$ is a square-free integer, and each prime divisor of $h_K$ is congruent to $2$ modulo $3$. Thus, Theorem \ref{PHT2} applies and gives $3 \in \textbf{R}_{K/\Q}$.

\section{Appendix (Another proof of Theorem \ref{Main1})} \label{secF}

In this appendix we give an alternate proof of Theorem \ref{Main1} as sketched by Jiuya Wang.
\begin{proof}
Let \(S\) denote the set of rational primes dividing the discriminant \(D_K\);
thus \(|S|=\omega(D_K)\), where \(\omega(D_K)\) denotes the number of distinct prime divisors of \(D_K\).
Since the extension \(H_\ell(K)/\Q\) is abelian, by the Kronecker--Weber theorem, there exists an integer $m$
such that
\[
H_\ell(K)\subseteq \Q(\zeta_m).
\]

As the extension \(H_\ell(K)/K\) is unramified , it follows that \(H_\ell(K)/\Q\) is unramified outside \(S\). Therefore, we see that
\[
m=\prod_{p\in S} p^{a_p}.
\]

Passing to Galois groups yields an injection
\[
\Gal{H_\ell(K)/\Q}\hookrightarrow (\Z/m\Z)^\times
 \cong \prod_{p\in S}(\Z/p^{a_p}\Z)^\times .
\]
For each prime power \(p^{a_p}\), the \(\ell\)-Sylow subgroup of
\((\Z/p^{a_p}\Z)^\times\) is a product of at most two cyclic \(\ell\)-groups.
Consequently there exists an absolute constant \(C_1\) (indeed $\le2$) such that
\[
\ell\text{-rank}\big((\Z/m\Z)^\times\big)\le C_1\,|S|.
\]
Hence
\[
\ell\text{-rank}\big(\Gal{H_\ell(K)/\Q}\big)\le C_1\,|S|.
\]
By Artin reciprocity, $\Gal{H_\ell(K)/K}$
is canonically isomorphic to the $\ell$-primary part of $C\ell(K)$. Thus, we have
\[
\ell\text{-rank}\big(C\ell(K)\big)\le C_1\,\omega(D_K).
\]
Therefore
\[
|C\ell(K)[\ell]|=\ell^{\ell\text{-rank}(C\ell(K))}
 \le \ell^{C_1\omega(D_K)}.
\]
Finally, using the classical bound
\[
\omega(D_K)\ll \frac{\log D_K}{\log\log D_K},
\]
we conclude that for every \(\varepsilon>0\),
\[
\ell^{C_1\omega(D_K)}\ll_{\varepsilon,\ell} D_K^{\varepsilon}.
\]
This proves the theorem.
\end{proof}

\end{document}